\documentclass[12pt]{article}
\usepackage{amsfonts,amsmath,amssymb,amsthm,amscd,epsfig,epic,eepic}
\usepackage{graphics}

\evensidemargin = \oddsidemargin
\textheight = \paperheight
\advance\textheight by -2in
\textwidth = \paperwidth
\advance\textwidth by -2in

\makeatletter
\gdef\thmhead@plain#1#2#3{%
  \thmname{#1}\thmnumber{\@ifnotempty{#1}{ }#2}%
  \thmnote{ {\mdseries#3}}}
\let\thmhead\thmhead@plain
\makeatother
\theoremstyle{plain}

\newtheorem*{thm1}{Theorem 1}

\newtheorem{thm}{Theorem}[section]
\newtheorem{lem}[thm]{Lemma}

\newtheorem{prop}[thm]{Proposition}
\newtheorem{crr}[thm]{Corollary}

\theoremstyle{definition}
\newtheorem{dfn}[thm]{Definition}

\newtheorem{rmk}[thm]{Remark}

\newtheorem{assu}{Assumption}
\newtheorem{step}{Step}
\newtheorem{step1}{Step}

\theoremstyle{remark}

\def\alinea#1{\hfill\break%
  \hbox to \parindent{\hss{\upshape{\bf #1)}}\enspace}\ignorespaces}
\def\bul{\hfill\break\hbox to\parindent{\hss$\bullet$\enspace}\ignorespaces}

\def\from{\colon}

\def\mod{\operatorname{mod}}

\def\C{\mathbf{C}}

\def\D{\mathbf{D}}

\def\N{\mathbf{N}}
\def\Q{\mathbf{Q}}
\def\R{\mathbf{R}}
\def\S{\mathbf{S}}

\def\Z{\mathbf{Z}}

\def\Cap_#1{\bigcap\limits_{#1}}
\def\Cup_#1{\bigcup\limits_{#1}}
\def\ol{\overline}

\def\wt{\widetilde}

\def\p{\mathbf 0}
\def\bi{B(\infty)}

\makeatletter
\def\cqfdsymb{\relax\protect\ifmmode\else\unskip\nobreak\fi
\quad\hfill$\bgroup
\vcenter{\hrule\hbox{\vrule\@height.6em\kern.6em\vrule}\hrule}\egroup$}
\def\cqfd{\cqfdsymb\endtrivlist}
\gdef\rom#1{\leavevmode\skip@\lastskip\unskip\/%
        \ifdim\skip@=\z@\else\hskip\skip@\fi{\normalshape#1}}
\makeatother
\title{Cubic polynomials with a parabolic point}
\author{{\sc P. Roesch}\footnote {IMT, Toulouse, France}}
\date{\today}
\begin{document}

\maketitle \maketitle
\begin{abstract}
We consider the family of cubic polynomials with a simple parabolic fixed point. We prove that the boundary of the immediate 
 basin of attraction of the  parabolic  point
is  a Jordan curve and give a description of the dynamics.  \end{abstract}
\section{Introduction and basic properties}
In this article, we consider polynomials  of degree three with a parabolic fixed point of multiplier 
$1$ and focus on the topological properties of the Julia set.  We prove

\begin{thm1}\label{th:borddeB}
The boundary of the immediate  basin of attraction of the parabolic fixed point is a Jordan curve, or two in the particular case of  the polynomial $z^3+z$.\end{thm1}
This Theorem  is the generalization to the parabolic case of~\cite{Ro1} (see also the appendix of \cite{Ro2}). It is also a first step for understanding the parameter space of cubic polynomials with such a parabolic fixed point (see~\cite{Ro3}).

In section~\ref{s:kiwi} (Proposition~\ref{p:kiwi}), we prove the local connectivity at the parabolic fixed point 
(and at all its iterated inverse images) 
reformulating  an argument of J.~Kiwi (see~\cite{Kiwi}). Theorem~\ref{th:borddeB} is proved  in section~\ref{s:puzzles}. We use Yoccoz' result for puzzles jigsawed by graphs containing  ``periodic" accesses. These accesses are  constructed in section~\ref{s:acces}.
Finally, we describe the dynamics of the Julia set in section~\ref{s:dyn}, it is based on results of section~\ref{s:wakes}.


\vskip 1em
For backgrounds in dynamics we recommand  the articles~\cite{DH,Mi1,CG}.
We fixed now some notations\,:
Every polynomial of degree three is conjugated (by an affine map) to a polynomial of the form $\lambda z+az^2+z^3$ with $a\in \C$
and $\lambda\in \C$ denotes the multiplier of the fixed point $0$. We restrict ourself to the cubic polynomials of the form 
$z+az^2+z^3$, {\it i.e.} possessing a parabolic fixed point of multiplier $1$ at $\p$.
  Let $f$ be such a polynomial. Denote by $K(f)$ its filled-in Julia set\,: $K(f)= \{z \mid f^{n}(z) \hbox{ is bounded}\}$, and  by $J=J(f)$, the Julia set of $f$\,:  its boundary. The map $f$ has two critical points in $\C$ and one at infinity.
The basin of attraction of $\p$, which is by definition $\tilde B=\{z \mid
f^n(z)\rightarrow \p\}$, contains at least one critical point. More precisely this critical point is in $B$, the {\it immediate} 
basin of $\p$, {\it i.e.} the connected components of  $\tilde B$ 
containing $\p$ in their boundary. Here there are two cases, depending on wether $a=0$ or not.

If  the parameter $a$ is $0$,  the parabolic point is double since $f(z)=z^3+z$.  
There are two connected components in $B$ (with the symmetry $z\mapsto -z$), each one contains a critical point. Thus the boundary of $B$  is the union of two Jordan curves since $f$ is sub-hyperbolic (see \cite{TY}).
In the case where  $a\neq0$, there is only one attracting direction so that $B$ has only one connected component (it is contained in $B$).
 If both critical points belong to $B$, it  is totally invariant, so that  $B=\tilde B$ is a topological disk and $\partial B$ is a Jordan curve since  $f$ is sub-hyperbolic as before.  For the same reason,  if one critical point belongs to  $\tilde B\setminus B$ or  to the basin of infinity, $ \bi=\{z \mid f^n(z)\rightarrow \infty\}$, or is mapped to the fixed point $\p$,    the boundary of $B$ is a Jordan curve, $\tilde B$ is a union of Jordan disks  (see~\cite{TY}).
Therefore we make the following assumption\;:
\begin{assu} Only one critical point is  in $\tilde B$, there is no critical point  in $\bi\setminus\{\infty\}$ and the orbit of the critical points avoids the fixed point $\p$. \end{assu}
Denote by  $c_0$ the critical point contained in $B$ and by  $c$ the other critical point (of $\C$).
\begin{rmk}Under assumption~1, the Julia set is connected (see~\cite{DH}) and $a\neq 0$ by the previous description. 
\end{rmk}
Let $\phi_\infty$ be the B\"ottcher coordinate at infinity. It is  defined from $B(\infty)$ onto $\C\setminus \overline \D$ and satisfies $\phi_\infty\circ f(z)=(\phi_\infty(z))^3$. For $v>0$ denote by   $E_\infty(v)$ the equipotentiel of level~$v$ : $\phi_\infty^{-1}(\{e^{v+2\pi i \theta}, \; \theta \in \S^1\} )$, and denote by $R_\infty(\theta)$ 
the external ray  of angle $\theta\in \R/\Z$\,: 
$\phi_\infty^{-1}(\{re^{2\pi i \theta},\; r>1\})\cup\{\infty\}$. Recall that every  ray with rational angle lands at a parabolic or repelling (eventually) periodic point (see~\cite{DH,Mi1}). Conversely, since $J(f)$ is connected one has (see~\cite{Mi1,Pe1})\,: 

\begin{prop}{\rm[Yoccoz]}\label{p:yocc}
For every repelling or parabolic (eventually) periodic point of $f$  there exists a ray $R_\infty(\theta)$ landing at it. Moreover,  the angle $\theta$ is rational.
\end{prop}
In particular, there is an external ray  landing at the parabolic fixed point  $\p$. This  ray is fixed since  $f'(\p)=1$ (see~\cite{GM, Pe1}). Thus,
 it is either $R_\infty(0)$ or $R_\infty(1/2)$ and, up to changing the B\"ottcher coordinate $\phi_\infty$ by $i\phi_\infty$,   we can suppose that  $R_\infty(0)$ lands at the parabolic fixed point $\p$.

\vskip 0.5 em
 Let $\phi_0$ be the Fatou coordinate at the parabolic fixed point $\p$, defined on a topological disk  $U\subset B$ (an attracting petal) such that $\phi_0:U\to\{z\mid \Re e(z)>0\}$  satisfies   $\phi_0(f(z))=\phi_0(z)+1$ and $\phi_0(f(c_0))=1$, 
  $\partial U$ contains the points  $\p$ and $c_0$,   and is smooth outside $\p$.   
  One can extend  $\phi_0$ to a
semi-conjugacy $\ol{\phi_0} \from B \to \C$, defined up to the post-addition of a constant.
The ``lines'' $E_\p(k)=\ol\phi_0^{-1} (\{k+it, t\in \R\})$ will play the role in the basin $B$ of the equipotentials of the usual  attracting case.


\section{Local connectivity  at the parabolic point.}\label{s:kiwi}
We define a sequence of neighbourhoods of the parabolic fixed point $\p$ as follows.
\begin{dfn}\label{d:kiwi}
Let $S_0$ be the graph $E_{\p}(1)\cup R_\infty(0)\cup\{\p\}$ and 
 $S_1$ be the connected component of $f^{-1}(S_0)$ containing $c_0$. Denote by  $U_0$  the connected component of $\C \setminus S'_1$ containing  $\p$ in its boundary, where $S'_1=S_1\setminus R_\infty(0)$.
Then denote by   $U_n$  the connected
component of  $f^{-1}(U_{n-1})$ containing $\p$ in its boundary and by $X$ the intersection  $\Cap_{n\ge 0} \ol{U_n}$.

\end{dfn}

Note  that there is   only one connected component of $f^{-1}(U_{n-1})$
containing~$\p$ in its boundary (otherwise  $\p$ would be a critical point),  so the sequence  $(U_n)$ is well defined .

\begin{rmk}\label{r:covering}
By construction, the sets
 $\ol U_n\cap \partial B$ are connected neighbourhoods of $0$ in $\partial B$.
Morever, for $n\ge 1$,  $U_n \subset U_{n-1}$ and the map  $f \from U_n \to U_{n-1}$ is a (ramified) covering of degree  $d\le2$.
\end{rmk}
\proof The critical point $c$ is not  in $U_n$ so the degree $d$ is less or equal than $2$. \cqfd

\begin{prop}~\label{p:kiwi}
Either $X$ reduces to the fixed point $\{\p\}$  or the rays $R_\infty(0)$ and 
$R_\infty(1/2)$ converge to $\p$ and  $\gamma=R_\infty(0)\cup R_\infty(1/2)\cup\{\p\}$
 separates $X\setminus \{\p\}$ from $B$.
\end{prop}
\begin{crr}\label{c:lc0}
The boundary $\partial B$ is locally connected at the parabolic point $\p$ and at all its iterated inverse images.
\end{crr}
\proof If $X$ reduces to the fixed point $\p$, the sequence $(\ol U_n\cap\partial B)$ is 
  a basis of connected neighbourhoods of $\p$ in $\partial B$.
Consider now the second case of Proposition~\ref{p:kiwi}. Let  $V$  be the connected
component of $\C \setminus \gamma$ containing $B$ and take  $V_n=V\cap U_n$.
The sequence $(\ol{V_n})$ has diameter shrinking to zero
since $\Cap_n \ol V_n\subset\ol{V}\cap (\Cap_n\ol U_n)=\ol{V}\cap X=\{\p\}$. Thus 
$(\ol V_n\cap\partial B)$  is a basis of connected neigbourhoods of the fixed point 
$\p$ in $\partial B$ since  each $\ol{V_n}$ intersects $\partial B$ under a connected 
set (by construction).

Let $x\in\partial B$ be  an iterated  inverse image of $\p$, {\it i.e.} there exists    $i$  such that $f^i(x)=\p$. 
Let  $W_n$  be the connected components of  $f^{-i}(V_n)$  containing $x$ in its closure. Then 
 $(\ol W_n\cap \partial B)$ form a basis of  connected neighbourhoods of $x$ in $\partial B$.
\cqfd

{\it Proof of Proposition~\ref{p:kiwi}}\ : It goes by contradiction. 
Assume that $X$ is not a point. Therefore it is a non trivial full
connected compact set (the intersection of
a decreasing sequence of closed disks). Hence, we can consider a
conformal representation $\Phi$ from $\C \setminus X$ onto $\C
\setminus \ol{\D}$. 

%
\begin{figure}[!h]\vskip 0cm
\centerline{
\input{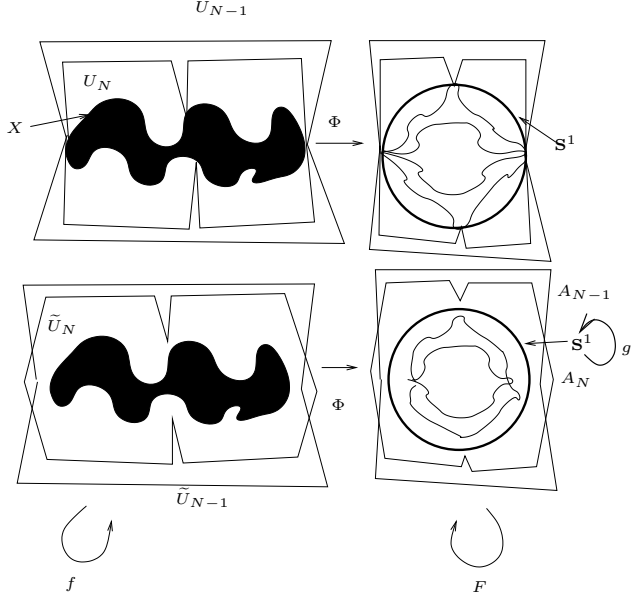}
} \caption{The neighborhood $U_N$ and a schematic representation of the map $g$ and of $X$
}\label{f:kiwi}
 \end{figure}

\begin{step}There exist  $A_N, A_{N-1}$,  neighbourhoods of
$\S^1$, such that $F=\Phi\circ f\circ \Phi^{-1}\from A_N \to A_{N-1}$
is  a holomorphic covering of degree $d\le 2$.
\end{step}
\proof Since the sequence $(U_n)$ is decreasing to $X$, there exists $N$ such that $c\notin U_N\setminus X$.  Therefore  $f^{-1}(X)\cap  U_N=X$ noting that  the set $X$ is fixed by $f$  (since $f(U_n)=U_{n-1}$ for all $n>0$).
Thus the map  $F=\Phi\circ f \circ \Phi^{-1}$ is  well defined from $\Phi(U_N\setminus X)$ 
to $\Phi(U_{N-1}\setminus X)$.  The open set   $\Phi(U_N\setminus X)$  is not a neighbourhood of $\S^1$ since $X\cap\partial U_N$ contains the iterated inverse 
images of the parabolic fixed point~$\p$.
Therefore, consider  $\wt U_{N-1}$  a small enlargement of $U_{N-1}$ near the points of  $X\cap \partial U_{N-1}$ (by a small ball for instance) such that $f(c)\notin \wt U_{N-1}\setminus X$ and 
 let   $\wt U_{N}$ be the connected component of $f^{-1}(\wt U_{N-1})$  containing~$\p$.
Thus, $f\from  \wt U_{N} \to \wt U_{N-1}$ is a (ramified) covering of degree $d\le2$, and if $W_i= \Phi(\wt U_i\setminus X)$, the map
$F\from W_{N}\to W_{N-1}$ is  then a  covering of degree $d$.
Now,  let $A_N$ be $W_N \cup \S^1 \cup \tau W_N$  where
 $\tau$ is the reflection through the circle\,; $F$ extends  to a  covering $F\from A_{N}\to A_{N-1}$  (of degree $d\le 2$) using Schwarz reflection principle . 
\begin{step} The restriction $g=F_{|_{\S^1}}$   is a degree $d$   covering preserving the
orientation, that can be lifted to an increasing homeomorphism  $\wt g\from\R\to\R$ with $\wt g(0)=0$.
Every  fixed point of $f$ in $X$ with a fixed access from $\C\setminus X$ 
gives a fixed point for $g$. 
\end{step}

\proof The map $F_{|_{\S^1}}$   is a degree $d$   covering preserving the
orientation since  $F\from A_{N}\to A_{N-1}$ is a holomorphic covering of degree $d$  preserving $\S^1$.  Let $\delta\subset\C\setminus X$  be a fixed access converging to 
 a fixed point $x\in X$. The access $\Phi(\delta)$ is fixed by $F$ and converges 
(see \cite{CG,Po}) to some point $p\in\S^1$ which is necessary fixed by $g$. 

\begin{step}\label{l:revtdeS1} Every fixed point $p$ of $g$ is weakly repelling, {\it i.e.} satisfies  for every $x\neq p$, 
$| g(x)-p|_{\S^1}>|x-p|_{\S^1}$.
Thus $g$ has exactly $d-1$ fixed points.
\end{step}
\proof  
Assume, to get a contradiction, 
that $p$ is attracting  on one side (at least). 
Then $p$ is  an attracting or a parabolic fixed point for  the holomorphic map $F$ so that there exists an arc 
$\alpha\subset W_N$ from $p$ to some point $q\in \S^1$  bounding an open set 
$\Omega_1$ in $\C\setminus \overline \D$ such that $F(\Omega_1)\subset \Omega_1$ (in the attracting domain of $p$).
Let $\Omega$ be $\Phi^{-1}(\Omega_1)$. It satisfies $f(\Omega)\subset \Omega$, thus the family $(f^n)$ 
is normal on $\Omega$. This gives morally  a ``neighbourhood'' of points of $J$ on which $f$ is normal. We prove now that this is impossible. There is a  fixed Fatou component $\wt\Omega\neq B(\infty)$  containing  $\Omega$  and (by Denjoy-Wolf's Theorem) $\ol{\wt\Omega}$ contains a fixed point attracting every point of $\wt\Omega$.
Moreover $\partial \Omega \cap X\subset\partial \wt \Omega$ contains  more than one point. So there is a cross cut $c$ of $X$ in $\tilde \Omega$. Indeed, the map $\Phi^{-1}$ admits  limit points at 
almost every $\theta \in(p,q)$  by Fatou's Theorem\,; if the limits
points were equal, say  to $0$ (to fixe the ideas), then the map 
$H(z)=\Phi^{-1}(z)\Phi^{-1}(\rho z)... \Phi^{-1}(\rho^k z)$ which is holomorphic would   admit $0$ as limit for almost every $z\in\S^1$
(were $\rho$ is a rotation such that $\Cup_{i\le k}\rho^i (p,q)$ cover all the circle)\,; thus $H$ would be  identically $0$ which contradicts the fact that $\Phi^{-1}$ is not constant.
This  cross cut  $c$  bounds a domain $U$.
Every point  of $\partial U\cap\partial \wt \Omega$ is in $X$ (since the Julia set of a polynomial is full). Hence $\partial \wt \Omega\subset X$ (by iterations of $U$).
This is not possible since $X$  is a full compact connected set  and $ \wt\Omega$ is a topological disk.

We prove now that $g$ has $d-1$ fixed points.
Since every  fixed point of $g$ is weakly repelling, the   graph of $\wt g$ has to cross the
diagonal and its translates from below to above. Therefore the
graph crosses exactly $d-1$ translates of the diagonal (including
the diagonal) exactly once. The result follows.
 \cqfd
\begin{step} The rays $R_\infty(0)$ and $R_\infty(1/2)$
converge to the fixed point~$\p$ and $\gamma$ separates $X\setminus\{\p\}$ from~$B$.
\end{step}
\proof By choice of the B\"ottcher map, the ray $R_\infty(0)$ converges to the parabolic fixed point~$\p$. It corresponds to him  a fixed point $x_0$ of $g$ in $\S^1$ (by Step~2).  This implies that  $d=2$ (Remark~\ref{r:covering}) since  by  Step~\ref{l:revtdeS1}  $d-1\ge 1$. Hence the map $g$ has exactly one fixed point.  We prove now that the landing point of $R_\infty(1/2)$ gives also a fixed point for  $g$ so we obtain  that $R_\infty(1/2)$  lands at $\p$. For this  we prove  that $R_\infty(1/2)$  enters every $U_n$ ($n>0$).

 Since $d=2$,   the critical point $c$ belongs to $ X$. Hence, there is a preimage  $R_\infty(\eta)$ of $R_\infty(0)$  ($\eta\in \{1/3,2/3\}$) converging to  $y$, a preimage of $\p$ in $ X$, which enters every  $U_n$ ($n>0$). The boundary of $U_n$ contains  external rays $R_\infty(\theta_n)$, $R_\infty(\theta'_n)$ landing at points of $\partial B$ and such that $\theta'_n<\eta<\theta_n$ (for the cyclic order on $\S^1$). Moreover,   these external rays satisfy $f(R_\infty(\theta_n))=R_\infty(\theta_{n-1})$, $f(R_\infty(\theta'_n))=R_\infty(\theta'_{n-1})$ and
$0<\theta_n<\theta_{n-1}<\theta'_{n-1}< \theta'_n<1$ (by definition of $U_n$).
Therefore, the sequences $(\theta_n)$ and $(\theta'_n)$
converge to limit angles $\theta$ and $\theta'$, which are distinct since  $\theta'\le\eta\le\theta$
and  since $\eta$ is prefixed while  $\theta$, $\theta'$  are fixed by multiplication by $3$.
Thus say  $\theta=0$ and $\theta'=1/2$. Therefore  the ray $R_\infty(1/2)$ enters every $U_n$. It  lands at a fixed point of $f$ in $X$, so gives  a fixed point on $\S^1$ for $g$ (Step~2).  But the map $g$ has  only one fixed point, by step~\ref{l:revtdeS1} (since $d=2$), so  $\Phi(R_\infty(0))$ and $\Phi(R_\infty(1/2))$ land at the same point.  Thus 
 $R_\infty(0)$ and  $R_\infty(1/2)$ both land at~$\p$ (see~\cite{Po}). 

We prove now that $\gamma$ separates $B$ from $X$. Let $V$ be the connected component of $\C\setminus \gamma$ containing $B$ and $W$ the other one. 
The critical point $c$ belongs to $W$ as well a the preimage of $B$, so
$X \cap W\neq\emptyset$. The sequence $V_n=U_n\cap V$ is decreasing. Since there is no   critical points in $V_0$, the map $f\from V_{1}\to V_{0}$  is an homeomorphism and  every point  of $V_0$ converges by iteration under $f^{-1}$
to a fixed point of $\overline{V}$ (by Denjoy Wolff's Theorem on each connected component of $V_0$). 
The map $f$ sends $X\cap V_1$ to  $X\cap V_0$, so the previous behaviour is not possible 
($X\cap V\subset f^{-n}(V_1)$) and  necessarily $X\cap V=\emptyset$.
\cqfd


\section{Access to the boundary of $B$.}~\label{s:acces}
The attracting petal  $\Omega=\phi_0^{-1} (\{x+it, x>1, t\in \R\})$  contains the critical value $f(c_0)$ and  satisfies $f(\Omega)\subset \Omega$.  We construct accesses to $\partial B$ lying in $B\setminus \Omega$.

\begin{dfn}\label{d:access}
A $k$-periodic access of angle $\theta$ is a  simple curve $\gamma(\theta)\subset B\setminus \Omega$ 
such that :
\begin{enumerate}
 \item$f^i(\gamma(\theta))\cap(B\setminus \Omega)$  are simple and  disjoint curves for $0\le i\le k-1$\,; 
\item $f^k(\gamma(\theta))\cap(B\setminus \overline  \Omega)=\gamma(\theta)$\,;
\item $f^i(\gamma(\theta))$  intersects $E_{\p}(1)$ at exactly one point
 denoted $y_i(\theta)$ for $0\le i\le k-1$\,; 
\item  the points $y_0(\theta), \ldots, y_{k-1}(\theta)$ are 
in the same cyclic order on the line $E_{\p}(1)$ than the points $z_j=e^{2i\pi \theta 2^j}$ 
for $0\le j\le k-1$. 
\end{enumerate}
 \end{dfn}
\noindent The $k$-periodic accesses   can be think of as   internal rays (in the case of an attracting basin). In particular, there exist periodic accesses for any angle $\theta$ periodic by multiplication by $2$. We construct them using their itinerary with respect to the following partition.  Recall that $S_1$ is the connected component of $f^{-1}(E_\p(0)\cup R_\infty(0)\cup\{\p\})$  containing $c_0$.
\begin{dfn} Let $\Xi_0$, resp. $\Xi_1$, be the unbounded connected component  of $\C\setminus S_1$
 which is on the left, resp. on the right,  of $R_\infty(0)$ when one goes from $\p$ to $\infty$.\end{dfn}

\begin{figure}[!h]
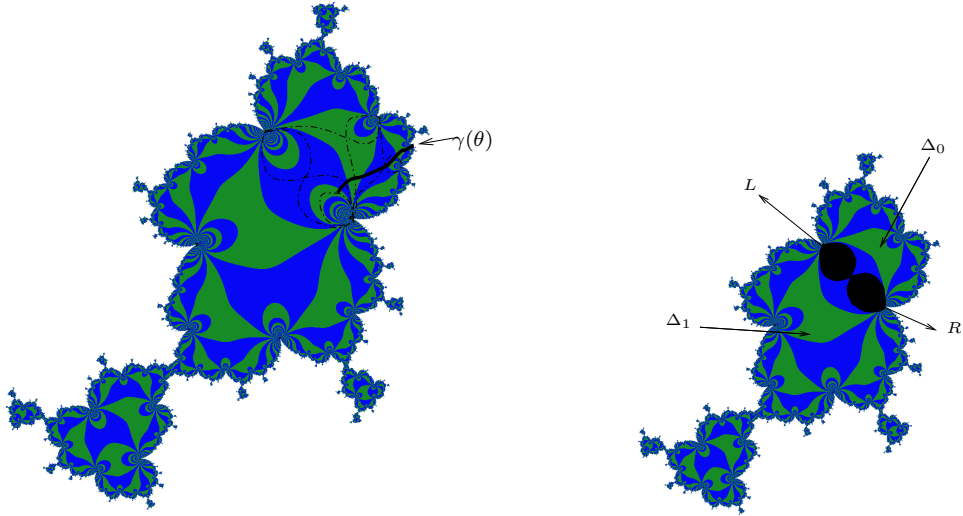
\vskip 0cm
\centerline{\input{gamma1.pstex_t}\hskip -2em \input{cbreduced.pstex_t}
} \caption{The construction of a $k-$periodic access.
}\label{f:}
 \end{figure}
 
For our purpose we do not need this construction  in the general case but only for angles of the form $\pm \frac{1}{2^m-1}$, so\,:

\begin{lem}\label{l:construction}
There exist $k$-periodic accesses of angle $\theta_\pm=\frac{\pm 1}{2^k-1}$. They  land at $k$-periodic points of $\partial B$.
\end{lem}
\proof  Let us do the construction fo $\theta_+$,  the argument is similar for $\theta_-$.

The set $f^{-1}(\Omega)$ is bounded by $E_{\p}(0)$, it has two connected components. Denote by $R$ the one containing the parabolic point $\p$ in its closure and by $L$ the other one. Let $\Delta_0$, resp. $\Delta_1$, be the connected component of $B\setminus E_{\p}(0)$ contained in $\Xi_0$, resp. $\Xi_1$.
The map $f\from B\setminus \overline{R\cup L}\to  B\setminus \ol\Omega$ admits two inverse branches $f_0\from  B\setminus \ol\Omega\to\Delta_0$
and $f_1\from B\setminus \ol\Omega\to\Delta_1$ which extend continuously  to the boundary by  $f_0\from  B\setminus \Omega\to\ol\Delta_0$
and $f_1\from B\setminus \Omega\to\ol\Delta_1$.  

Let  $y_0(\theta)$ be a point on $E_{\p}(1)$ with  $\Im m(\phi(y_0(\theta)))>0$ and denote by  $y'_0$ the point $g(y_0(\theta))$ where $g=f_0^{k-1}\circ f_1$ is defined on $B\setminus \Omega$.
Consider the set $U=\Delta\cup f_0(\Delta)\cup\cdots\cup f_0^{k-1}(\Delta)$ on which the extension of the Fatou coordinate $\overline \phi$ is an homeomorphism to $[-(k-1),1]\times \R^+$, with  $\Delta=\overline {\phi_0^{-1}(]0,1]\times \R^+)}$.
Let $\delta_0$ be the straight line joining $\phi_0(y_0(\theta))$ and $-(k-1)=\phi_0(f_0^{k-1}(c_0))$. Thus $\phi_0^{-1}(\delta_0)$ joins $y_0(\theta)$ with  $f_0^{k-1}(c_0)$,
 the iterates $f^j(\phi_0^{-1}(\delta_0))$ are disjoints and the points $y_i(\theta)=f^j(\phi_0^{-1}(\delta_0))\cap E(1)$ are, in the Fatou coordinate, the intersection of $\delta_0+j$ with the line $x=1$, so are ordered like the points $z_j$.
Let $\delta_1$ be any closed arc in $L$ joining $f_1(y_0(\theta))$ and $c_0$.
Consider the curve $\gamma(\theta) =\cup_{i\ge 0} g^i(\delta)$ where $\delta=\phi_0^{-1}(\delta_0)\cup f_0^{k-1}(\delta_1)$. 

By construction, $\gamma(\theta)$  is a curve  since $\ol \delta$ joins $y_0(\theta)$ and $g(y_0(\theta))$. It has no self-intersection. There are none on $\delta$ by choice of $\delta$ and any point of $\delta\cap \ol{g^i(\delta)}$ would give after composition by $f^i$ a point in $\delta\cap\ol\Omega$ which is empty.
For the same reason, since the $f^i(\ol \delta)$ are disjoint, the iterates  $f^i(\ol \gamma(\theta))\cap B\setminus \Omega$ are disjoint.
Finally, the cyclic order on $E_{\p}(1)$ is given by the order on the line $\{x=1\}$ in the Fatou coordinate. Thus $\gamma(\theta)$ is a $k$-periodic access of angle $\theta$. 

There exists some point $x\in\partial B$ such that the iterates of $g\from\Omega\to\Omega$ converge uniformly to $\{x\}$ on every compact set of $\Omega$ (by Denjoy Wolff Theorem). Applied to $C=g(\delta)$ it follows that $\gamma(\theta)$ converges to $x$.
By construction $f^i(\gamma(\theta))\subset \Xi_0$ for $0\le i\le k-2$, $f^{k-1}(\gamma)\subset \Xi_1$ and $f^k(\gamma(\theta))\subset \gamma\cup\ol\Omega $ so $x$ is exactly $k$-periodic.
\cqfd


\section{Wakes  of $B$.}~\label{s:wakes}
The construction of the wakes of a basin is classical (see~\cite{Mi2}). It allows to localize 
points with respect to external rays  or internal accesses.
\begin{lem}\label{s:semiconj}
There exists a map $\Theta\from \partial B \to\R/\Z$ such that $\exp(2 i \pi \Theta)\from \partial B \to \S^1$ is a  semi-conjugacy between $f$ and  $z\mapsto z^2$. Moreover
$\Theta(z)$ is dyadic if and only if $z$ is an iterated preimage of $\{\p\}$.
\end{lem}
\proof
Let $S_1$ be the connected component of $f^{-1}(E_{\p}(1)\cup \{\p\}\cup R_\infty(0))$ containing $c_0$. Let $\Xi_0$, resp. $\Xi_1$, be the connected component $\C\setminus  S_1$ containing $\Delta_0$, resp. $\Delta_1$.
For $z\in L$,  let $\Theta(z)$ be   the sum $\sum_{k\ge 0}\frac{\epsilon_k}{2^k}$ modulo $1$ where $\epsilon_k=0$, resp. $1$, if $f^k(z)\in \Xi_0$, resp $\Xi_1$. Remark that $\Theta(z)$ is not  dyadic otherwise some $f^i(z)$ would be in $\partial B\cap\Cap_n U_n=\{\p\}$ (where $U_n$ are defined in Definition~\ref{d:kiwi}).
 Now, set $\Theta(\p)=0\mod 1$ and  $\Theta(\beta)=1/2\mod 1$ where $\beta$ is the first preimage of $\p$ (on $B$). For any  iterated inverse image  $z$ of~$\beta$,  with   $f^k(z)=\beta$, the sequence $\epsilon_k$  depends on the position of $f^i(z)$ relatively to $\Xi_0$ and $\Xi_1$ for $i<k$ and continues with the sequence $1\ol 0$ or~$\ol 1$.
Thus, the map  $\exp(2i\pi \Theta)$ is well defined on $\partial B$ and  semi-conjugates $f$ to $z\mapsto z^2$.
 \cqfd
\begin{lem}\label{l:limbs} The map $\Theta$ extends to $\Theta\from  \C\setminus B \to\R/\Z$  such that $\exp(2 i \pi \Theta)\from \C\setminus  B \to \S^1$ is a   semi-conjugacy between $f$ and  $z\mapsto z^2$. This map  is constant along external rays and throughout any region of $\C\setminus B$ which is bounded by two external rays landing at a common point. 
\end{lem}
\proof
Let $S_n$ be the connected component, containing $c_0$, of $f^{-n}(E_{\p}(1)\cup \{\p\}\cup R_\infty(0))\setminus \cup_{i\le n} f^{-i}(c)$ (if the orbit of $c$ intersects $S_1$). Any point  $z$ of $\C\setminus \overline B$,  which is not on an external ray $R_\infty(\theta)$ with $\theta$ dyadic, is contained in  a  connected component $T_n(z)$
of $\C\setminus S_n$. 
These decreasing  components $T_n(z)$ intersect $B$ so that  there exists a point $u\in \partial B$ that belongs to  $\ol{T_n(z)}$ for every $n\ge 0$. Define $\Theta(z)=\Theta(u)$.
If $z$ belongs to some external ray   $R_\infty(\theta)$ with $\theta$ dyadic, define 
$\Theta(z)=\Theta(u)$ where $u$ is the landing
point of $R_\infty(\theta)$ (iterated inverse image of~$\p$).
 \cqfd

 \begin{dfn}Denote by  $W(z)$ the  {\it wake} of $z\in \C\setminus  B$\,:  the connected component of $\Theta^{-1}(\Theta(z))$ containing $z$. 
\end{dfn}

\begin{lem}\label{l:unicity}For any $z\in \C\setminus B$ there is a $u\in \partial B$ and 
unique $\theta$ such that $z$ belongs to $W(u)$ with $\Theta(z)=u$. 
Thus we can write without ambiguities $W(z)=W(\theta)$.
\end{lem}
\proof  In the proof of Lemma~\ref{l:limbs} $W(z)=W(u)$ for some $u\in \partial B$. Assume, to get a contradiction,  for some $v\in \partial B$ such that $W(v)=W(u)$, with $\Theta(u)=\theta$ and $\Theta(v)=t$. By definition  of $\Theta$, these points 
are in different connected components of $\C\setminus  S_n$ since there is a dyadic angle separating $\theta$ and $t$. This leads to a contradiction.
\cqfd

\begin{rmk}\label{r:imwake} The image of a wake $W(t)$ not containing the critical point $c$ in its closure is the  wake $W(2t)$.\end{rmk}

\begin{lem} If $c\in \ol {W(t)}$ then the image $f(W(t))$ is the whole plane. Moreover $f(c)$ belongs to $\ol {W(2t)}$.
\end{lem}
\proof  If the closure of the  wake contains  the critical point, it also contains the preimage $B'$ of $B$ and the preimage of all the wakes attached to $B$. Clearly its image covers all $\C$.
Let $t_n,t'_n$ be the angles of the rays bounding $T_n(c)$. The wake $W(c)$ contains $B'$ the 
inverse image of $B$ so  contains also a  preimage of these two external rays. 
Let $V_n$ be the connected component containing $c$ bounded by these four external rays 
together with $\ol B$  and  $\ol B'$.  It is mapped with degree 
two on the component of $\C\setminus (\ol B\cup R_\infty(2t_n)\cup R_\infty(2t'_n))$ containing $f(c)$. 
The result follows.
\cqfd
\begin{lem}\label{l:wakeimage}
For any wake $W$ not reduced to the closure of an external ray, there exists a smallest $n\ge 0$ such that the  closure of $f^n(W)$ contains the critical point $c$.\end{lem}
\proof The size of a wake $W(z)$ can be defined as the limit of the difference of the angles of the external rays  of $T_n(z)$.
If a  wake does not contain the critical point  in its closure, its size is multiplied by $3$ at each iteration.  Hence we deduce from Remark~\ref{r:imwake} that  the closure of $f^i(W)$ contains the critical point $c$ for  some $i\ge 0$.   
\cqfd

\begin{lem}\label{l:critperiod}
If $c\notin \partial B$,  there exists some $t\in \Q/\Z$, $k$-periodic by multiplication by $2$, such that $W(t)$ contains $c$.\end{lem}
\proof Let $W(t)$ denotes the wake containing the critical point. Since $c$ belongs to $K(f)$ it is not reduced to the closure of an external ray.
The critical  value $f(c)$ belongs to the wake $W=W(2t)$. If $t\neq 0$, $W$ does not contain the critical point $c$ so there is some $n\ge 1$ such that $f^n(W)$ contains  $c$ in its closure (Lemma~\ref{l:wakeimage}).
Therefore $f^n(W)=W(t)$ and since $f^n(W)$ is a wake of angle $2^n(2t)$, we get by Lemma~\ref{l:unicity} that $2^{n+1}t=t\mod 1$.
\cqfd

\section{Puzzles.}~\label{s:puzzles}
We can now use Yoccoz result with graphs containing the $k$-periodic access constructed
 in section~\ref{s:acces}.
  Let $x(\theta)$ be the landing point of $\gamma(\theta)$ and  take  an external ray  $R_\infty(\zeta)$ converging to $x(\theta)$ (Proposition~\ref{p:yocc}). The ray $R_\infty(\zeta)$ has period $k$ exactly since the point $x(\theta)$  has rotation number $1$ (see~\cite{GM,Pe1}). Thus $\Cup_{0\le i\le k-1}f^i(\gamma(\theta)\cup R_\infty(\zeta)\cup x(\theta))$ is forward invariant in $\C\setminus \ol \Omega$.

\begin{dfn}\label{d:graph}
Let $\Gamma(\theta)$ be the following  graph  $$\Gamma(\theta)=E_\infty(1)\cup E_{\p}(1)\cup \Cup_{0\le i\le k-1}
f^i(\gamma(\theta)\cup R_\infty(\zeta)\cup x(\theta))\cap (B\setminus \Omega).$$
\end{dfn}

\begin{dfn} Let $\Gamma_n=f^{-n}(\Gamma(\theta))\cap B\setminus \Omega$. The puzzle pieces of depth $n$ are  the connected component of $\C\setminus \Gamma_n$\,;  denote by $P_n(z)$ the one  containing $z$.
\end{dfn}

\begin{prop}\label{p:condition}
Let $L$ be an invariant subset of $J(f)$ such that for every point  $z$ of $L \cup \{c\}$,  there exist a graph $\Gamma$  and a sequence $(n_i)$ such that  the puzzle pieces $(P_{n_i}(z))_{i\in\N}$, defined by the graph  $\Gamma$,  satisfies   $\ol{ P_{n_i}(z)}\subset P_{n_i-1}(z)$. Then for every point $z\in L$ either $\cap P_n(z)=\{z\}$ or there exists $j>0$ such that $f^j\from P_{n+j}(c)\to P_n(c)$ is renormalizable and  $\cap P_n(z)$ is {\rm(}an iterated inverse image of{\rm)} the  filled-in Julia set $K(f^j)=\cap P_n(c)$ of  the renormalization.
\end{prop}
Recall the definition of a renormalizable map in our context\,:
\begin{dfn}The map $f$ is $k$-renormalizable, or $f^k$ is renormalizable, near the critical point $c$,  if there exist
two open  topological disks $U,V$  with $U$ compactly contained in $V$ and  such that $f^k\from U\to V$ is a ramified double covering with $f^{kn}(c)\in U$ for every $n\ge 0$. Then,  
the filled-in Julia set of the renormalization  is $K(f^k):=\cap_{n\ge 0}f^{-kn}(U)$. 
\end{dfn}
{\it Proof of Proposition~\ref{p:condition}.} The Proposition follows from the work of Yoccoz.  We refer here to the exposition of~\cite{Ro1}. Recall that $z$ is said to be {\it non recurrent} to $c$ if $c\notin P_n(f^i(z))$ for every $i\ge 0$ and for every $n\ge n_0$ for some integer $n_0$.  Lemma~1.22 and Remark~1.23 solve the case of non recurrent points where we do not need to have a non degenerate annulus around the critical point.
In the case where $z$ is recurrent to $c$, we only need to have a non degenerate annulus surrounding the critical point $c$, but of arbitrary large depth.  Lemma~1.25 deals with the case of non persistently recurrent points, two cases appear. Here we have an infinite sequence  $(n_i)$ such that  $\ol{ P_{n_i}(c)}\subset P_{n_i-1}(c)$ so we fall in the first case, {\it i.e.} $p\ge l$.
Finally, in Lemma~1.27, where the persistently recurrent points are treated, the actual condition is also enough since we pull-back the infinite sequence $P_{n_i}(c)$ to puzzle pieces around $z$.\cqfd
\begin{figure}[!h]\vskip 0cm
\centerline{
\input{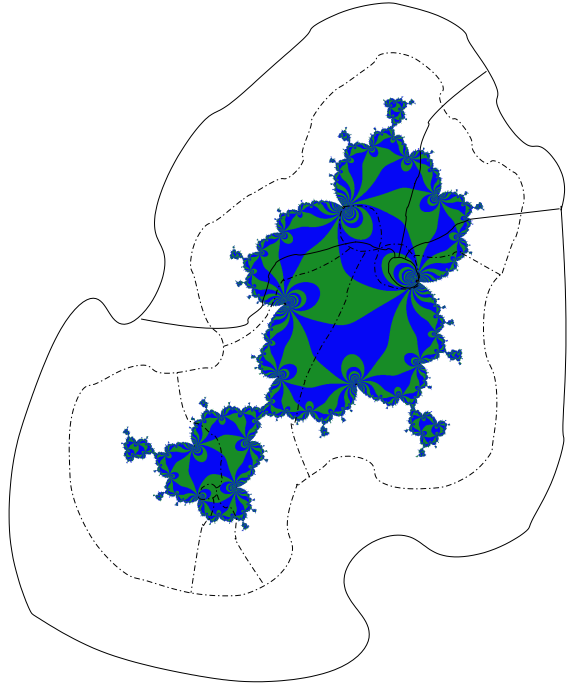}
} \caption{A graphe $\Gamma$ and its inverse image $\Gamma_1$
}\label{f:graphe}
 \end{figure}

\begin{lem}\label{l:condition} If $c\notin W(\p)$, the hypothesis of Proposition~\ref{p:condition} are  satisfied for the set $\displaystyle{L= \partial B\setminus   \cup_{i\ge 0} f^{-i}(\p)}$,  by graphs $\Gamma(\theta_\pm)$ with  $\theta_\pm=\frac{\pm 1}{2^k-1}$ and $k\ge 0$.
\end{lem}
\proof The set $L$ is clearly invariant by $f$. Step 1, 2, 3 do not use the assumption
$c\notin W(\p)$.

\begin{step1}\label{s:anneau} We consider the angle $\theta=\theta_+$. Let $P_0^+$ be the puzzle piece containing $\p$ in  its closure. Let $P^+_1$ be one of the puzzle pieces of depth $1$ intersecting $\Delta_1$ and not containing $\p\cup f^{k-1}(\gamma(\theta))$ in its closure.
 The annulus $P^+_0\setminus \ol P^+_1$ is non degenerate.
\end{step1}
\proof
Denote by $T$ the cycle of accesses $\{f^i(\gamma(\theta))\cap (B\setminus \Omega), \ i\in\{1,\cdots,k-1\} \}$. By construction, the puzzle piece $P^+_1$ is bounded in $B$ by two accesses $\gamma_i,\gamma_j$ of $f^{-1}(T)\setminus (T\cup\gamma(\theta))$.
and a portion of $E_{\p}(0)$ (see figure~\ref{f:graphe}). The curves $\gamma_i,\gamma_j$
 do not intersect $\partial P_0^+$ since on one hand, any intersection with an element of $T$ would give (after iterations) an intersection between  $\gamma(\theta)$ and  some element  of $T$\,; and on the other hand, $\gamma_i$, $\gamma_j$  start  on $E_{\p}(0)$ and  stay in $B\setminus (L\cup R)$ so cannot cut $E_{\p}(1)$. Remark that $E_{\p}(0)$ intersects $E_{\p}(1)$ at $\{\p\}$ but this point is not in the closure of $P^+_1$. Now, the portion of $E_{\p}(0)$ in $\partial P_1^+$ do not intersect $\partial P_0^+$. Indeed, the points $y'_{k-1}, \gamma_i\cap E_{\p}(0), \gamma_{j}\cap E_{\p}(0)$ are on $\partial L$ in the same cyclic order than $y_0(\theta), y_i(\theta),y_j(\theta)$, so the part of $E_{\p}(0)$ on the boundary of $P^+_1$ is included in 
 $\partial \Delta_1\cap  E_{\p}(0)$  after $y'_{k-1}$  so cannot intersect  $\Gamma$.
Thus $\ol P_1^+\subset P_0^+$. \vskip 0.5em
 The analogue  statement for  $\theta_-$ insures that the annulus  $P^-_0\setminus \ol P^-_1$ is non degenerate.

\begin{step1}\label{s:landingpt}
The landing point $x_j(\theta)$ of the access $f^j(\gamma(\theta))$ satisfies $\Theta(x_j(\theta))= 2^j \theta \mod 1$.
\end{step1}
\proof
 The map  $\exp(2i\pi\Theta)$  semi-conjugates  $f$ to  $z\mapsto z^2$,  so that $\exp(2i\pi\Theta(f^k(x_0(\theta))))=\exp({2^k}2i\pi\Theta(x_0(\theta)))$.  Since $x_0(\theta)$ is $k$-periodic, $f^k(x_0(\theta))=x_0(\theta)$, so that  
$\Theta(x_0(\theta))={2^k}\Theta(x_0(\theta))\mod 1$ and  $\Theta(x_0(\theta))$ is of the form $\frac{l}{2^k-1}$.
The points $x_0(\theta)$, $\cdots$ ,  $f^{k-2}(x_0(\theta))$ belongs to $\Xi_0$ and $f^{k-1}(x_0(\theta))$ is in $\Xi_1$, so it is easy to see that $\Theta(x_0(\theta))=\frac{1}{2^k-1}$ and the result follows.

\begin{step1}\label{s:infinite} For any point $z\in L$, there exist $k>0$ and $\theta\in \{\theta_\pm=\frac{\pm1}{2^k-1}\}$, such that the orbit of $z$ enters infinitely many times in 
one of the puzzle pieces, $P^\pm_1$,  of depth $1$  considered in Step~\ref{s:anneau}
for the graph $\Gamma(\theta_\pm)$.
\end{step1}
\proof  If $z$ is periodic we chose $k>0$  large  enough so that period of $\theta_\pm$ is larger than the one of $z$, so  $z\notin \Gamma(\theta_\pm)$. 
We use  the map $\Theta$ to localize the orbit of $z$.
Let $x'_j(\theta_\pm)$ be the preimage of $x_{j+1}(\theta_\pm)$ different from $x_j(\theta_\pm)$ for $0\le j\le k-2$.
By Lemma~\ref{s:semiconj} and Step~\ref{s:landingpt}, it is easy to see that 
$\Theta(x'_j(\theta_\pm))=2^j\theta_\pm+1/2\mod 1$. Hence, 
\begin{eqnarray*}\Theta(x_{k-1}(\theta_+))<\Theta(x'_0(\theta_+))=\theta_++\frac{1}{2}\mod 1<\Theta(x'_1(\theta_+))<\cdots<\Theta(x'_{k-2}(\theta_+))=\frac{3}{4}+\frac{\theta_+}{2}<1&&(*)\\
0<\Theta(x'_{k-2}(\theta_-))=\frac{1}{4}+\frac{\theta_-}{2}\mod 1<\cdots<\Theta(x'_0(\theta_-))=\theta_-+\frac{1}{2}\mod 1<\Theta(x_{k-1}(\theta_-)).&&(**)
\end{eqnarray*}
Assume that there exists a subsequence $(n_j)$ such that $2^{n_j}\Theta(z)\mod 1$ tends to $1/2$, then $2^{(n_{j}-1)}\Theta(z)\mod 1$  admits a sub-subsequence that tends to $1/4$ or to $3/4$.
Then the corresponding  sub-subsequence  of points $f^{(n_{j_l}-1)}(z)$  belongs in the first case (resp. in the second case) to $P_1^-$ (resp.   $P_1^+$), the puzzle piece containing $f^{-1}(\beta)\cap\Xi_0$ (resp. $f^{-1}(\beta)\cap\Xi_1$) by the inequality $(**)$, (resp. $(*)$), with $\beta=(f^{-1}(\p)\setminus\{\p\})\cap \partial B$.

Assume now that there exists $\epsilon>0$ such that for every $n\ge 0$ the angles 
$ 2^{n}\Theta(z)\mod 1$ avoid the interval $(1/2-\epsilon,1/2+\epsilon)$ containing $1/2$. Since $\Theta(z)$ is not dyadic,  there exists a subsequence $(n_j)$ such that $2^{n_j}t\mod 1$ belongs to the interval  $(1/2,3/4)$ (there exist infinitely many sequences of $10$ in its dyadic expansion). 
For $k$ large enough, $\theta_+=\frac{1}{2^k-1}<\epsilon$ so $\Theta(f^{n_j}(z))\mod 1$ belongs to the interval $(1/2+\theta_+,3/4)$. Therefore, by  the inequality  $(*)$,  the sequence $f^{n_j}(z)$ belongs to one of the  puzzle pieces $P_1^+$ considered in Step~\ref{s:anneau}.

\begin{step1}
Since $c\notin W(0)$, there exist $k>0$ and $\theta\in \{\theta_\pm=\frac{\pm1}{2^k-1}\}$, such that the orbit of $c$ enters infinitely many times in 
$P^\pm_1$,  one of the puzzle pieces of depth $1$  considered in Step~\ref{s:anneau}
for the graph $\Gamma(\theta_\pm)$.
\end{step1}
\proof Consider first the case where  the critical point $c$  belongs to  $\partial B$ in a wake $W(t)$.
If  the angle $t$ is  dyadic,  some iterate $f^i(c)$   is in $\partial B\cap W(0)$,  so coincides with $\{\p\}$ by section~\ref{s:kiwi}. This contradicts our assumption~1. If the angle $t$ is not dyadic, the result follows from Step~\ref{s:infinite}.

Consider now the case where the critical point is not on $\partial B$. By Lemma~\ref{l:critperiod} it  belongs to a wake $W(t)$ with $t$ periodic of period say $k$ by multiplication by $2$. Let  $u$ be a point of $W(t)\cap\partial B$. The points $f(u)$ and $f(c)$ belong to the  wake $W(2t)$. Therefore the points  $f^n(u)$ and $f^n(c)$  belongs to the same puzzle pieces of depth $0$ and $1$ for $n<k$. The result then follows from Step~\ref{s:infinite}.

\noindent{\bf End of the proof.}  Let $(n_j)$ be one of the sequences obtain in  Step~\ref{s:infinite}  such that $f^{n_j}(z)$ belongs to $P_1^\pm$. The map $f^{n_j-1}\from P_{n_j-1}^\pm(z)\to P_0^\pm$ is a (ramified) covering   which maps $P_{n_j}^\pm(z)$ to $P_1^\pm$.  Hence, the annulus 
$P_{n_j-1}^\pm(z)\setminus\ol{P_{n_j}^\pm(z)}$ is  a non degenerate  and therefore the assumptions of Proposition~\ref{p:condition} are satisfied.
\cqfd

\begin{lem}\label{l:w0}
If $c\in W(0)$,  for every point $z\in L$,  there exist $k>0$ and $\theta\in \{\theta_\pm=\frac{\pm1}{2^k-1}\}$, such that the sequence of  puzzle pieces $P_n(z)$ defined by this graph shrink to $\{z\}$.
\end{lem}
\proof Let $z$ be a point of $L$. By Step~\ref{s:infinite}, there exist $k>0$ and $\theta\in \{\theta_\pm=\frac{\pm1}{2^k-1}\}$  such that the orbit of $z$ enters infinitely many times in $P^\pm_1$,  one of the puzzle pieces of depth $1$  considered in Step~\ref{s:anneau}
for the graph $\Gamma(\theta_\pm)$.

If the orbit of $c$ is totally included in $W(0)$, the diameter of the iterated inverse images of  $P^\pm_1$ shrinks to $0$ by the so called "Shrinking Lemma"  (since 
$\ol  P^\pm_1$ do not intersect the post-critical set). So the result follows.

If $f^j(c)\in W(0)$ for $0\le j< k$ and $f^k(c)\in W(t)$ for some $k>0$ and  $t\neq 0$, then $f^{k-1}(c)$ belongs to the inverse image of $W(t)$ attached to $B'$ (the inverse image of $B$).
Therefore, the puzzle piece $P_1(f^{k-1}(c))$ is compactly contained in $P_0(c)$.
Indeed, it contains only internal accesses in $B'$ and do not intersect $B$\,; for the external rays and equipotential it is clear that there are no intersections.
Thus, $P_{k-1}(c)\setminus \ol P_k(c)$ is a non degenerate annulus surrounding $c$ (by pull back of the puzzle pieces).
Then the result follows from the alternative formulation of Proposition~\ref{p:yocc} which is the following.
For a polynomial $f$, if the critical point is surrounded by a non degenerate annulus, then for every point $z$ that is surrounded by infinitely many non degenerated annuli, the intersection $P_n(z)$ reduces to $\{z\}$ or to an iterated inverse image of the filled in Julia set a a renormalization of $f$ (as in Proposition)~\ref{p:yocc}. This formulation is exactly Theorem~1.10 of~\cite{Ro1}. If the map is renormalizable, its filled-in Julia set is contained in $W(0)$ (that contains $c$). Hence if the intersection of the puzzle pieces 
$P_n(z)$ is an inverse iterated image of this filled-in Julia set, the point $z$ would belong to a dyadic wake, because $c\in W(0)$. This is not possible for points in $L$.  \cqfd

Now we can prove the Theorem~\ref{th:borddeB} summarized in 
 \begin{thm1}The boundary $\partial B$ is locally connected. It is a Jordan curve.
\end{thm1}

\proof The boundary $\partial B$ is locally connected at every iterated inverse image of the parabolic fixed point  $\p$ by 
Corollary~\ref{c:lc0}. For $z\in L$,  Lemma~\ref{l:condition}, 
Proposition~\ref{p:condition} and Lemma~\ref{l:w0} imply that the intersection 
$\cap P_n(z)$ is equal  either to $\{z\}$ or to an iterated  inverse image (say by some $f^j$) of the filled in Julia set of some renormalization $f^l \from P_{n+l}(c)\to P_n(c)$. Note that the intersection $\partial P_n(z)\cap\partial B$ is clearly a connected set. Therefore,   $\partial B$ is locally connected at  $z$ in the first  case.

In the second case, the  puzzle pieces $P_n(c)$ intersect $B$ since  $f^j(z)\in P_n(c)$  for every $n\ge 0$ (and  $z\in \partial B$).  Let  $R_\infty(\zeta_n)$, $R_\infty(\zeta'_n)$ be  the two external rays of $\partial P_n(c)$ that land on $\partial B$. 
Since $P_{n+1}(c)\subset P_n(c)$, the angles
$(\zeta_n)$ and $(\zeta'_n)$ form  monotone and bounded sequences, they converge to limit angles $\zeta$ and $\zeta'$ respectively. Moreover, the periodicity property $f^l (P_{n+l}(c))=P_n(c)$ implies that the limits  $\zeta$ and $\zeta'$ are periodic by multiplication by $3$. Hence the rays $R_\infty(\zeta)$ and $R_\infty(\zeta')$ are fixed by $f^l$  so converge to  fixed points, with rotation number $1$. They belongs to $K(f^l)=\cap P_n(c)$ since $R_\infty(\zeta)$ and $R_\infty(\zeta')$ enter every  puzzle piece $P_n(c)$ for $n\ge 0$. There is only one fixed point of rotation number $1$ for the renormalized map $f^l$ in $K(f^l)$ (since $f^l$ is conjugated to a quadratic polynomial), denoted  $\beta(f^l)$. Therefore the two rays $R_\infty(\zeta)$ and $R_\infty(\zeta')$ both converge to $\beta(f^l)$.

 Let $V_n$, resp. $W_n$, be the connected components of 
$P_n(c)\setminus( {R_\infty(\zeta)\cup R_\infty(\zeta')}\cup\{\beta(f^l)\})$ intersecting $B$, resp. not intersecting $B$. We show  that  $(\ol V_n)$ forms a basis of neighbourhoods  of $\beta(f^l)=\partial B\cap K(f^l)$ in $\partial B$. For this we prove that $V_n\cap K(f^l)=\emptyset$. We only need to consider the case where 
$\beta(f^l)\neq \p$, because of  Proposition~\ref{p:kiwi}. In this case, note that   $V_n$ is connected. 
and that the wake $W_n$ contains the critical point $c$ (otherwise it would be mapped by $f^l$ homeomorphically to itself which contradicts Lemma~\ref{l:wakeimage}).  Hence, the map  $f^{l}$  admits an inverse branch on $V_n$, noted $g\from V_n\to V_{n-l}$.  There is some fixed point of $\partial V_n$ such that every orbit under $g^n$ converges to it  (by Denjoy-Wolff's Theorem). This implies that no point of  $K(f^l)$ is in $V_n$. Hence, the curve $R_\infty(\zeta)\cup R_\infty(\zeta')\cup \{\beta(f^l)\}$ separates $B$ from $K(f^l)\setminus\{\beta(f^l)\}$.
Therefore $\beta(f^l)=f^j(z)$. Let $Q_n$ be the pull-back of $V_n$ by the coresponding brach of $f^j$. The sequence $(\ol Q_n)$ forms  a basis of neighbourhoods of $z$ in $\partial B$.

In both cases we get that $N_n=\partial P_n(z)\cap\partial B$ form a basis of connected neighborhoods of $z$  in $\partial B$ (since $\cap N_n= \{z\}$). Thus, $\partial B$  is locally connected.
The classical theorem of Caratheodory  implies that $\partial B$ is a curve (since $B$ is a topological disk). Finally, this curve  is simple since  $K(f)$ is a  full compact set (by the maximum principle).
\cqfd

\begin{crr}\label{c:theta} The map $\exp(2i\pi\Theta)$ defined in Step~\ref{s:semiconj} of the proof of Lemma~\ref{l:condition} is a conjugacy between $f$ on $\partial B$ and $z\mapsto z^2$ on $\S^1$.
\end{crr}
\proof The map  $\exp(2i\pi\Theta)$ is a semi-conjugacy between $f$ and $z\mapsto z^2$. Two points with the same image have the same itinerary so are in the same iterated inverse image of $\Xi_0$ and $\Xi_1$ and thus in the same puzzle piece. Since they are on $\partial B$ they coincide.\cqfd

\section{Global description.}~\label{s:dyn}
\begin{lem}\label{l:2ray}For any point $z$ of $\partial B$ there are at most two external rays converging to it.\end{lem}\proof
Assume, to get a contradiction, that there are at least three rays landing at $z$. They  define two 
adjacent wakes $W_1,W_2$. Lemma~\ref{l:wakeimage} insures that either $f^n(z)=c$ for some $n$ or there exist $n_1,n_2$ such that $f^{n_1}(W_1)=f^{n_2}(W_2)$ is a wake containing the critical point. 
The second alternative is not possible, because of the cyclic order of the rays at $z$.
In the first alternative, $f$ identifies at most two rays, so that  at least two rays still land at $f(c)$. We can apply  Lemma~\ref{l:wakeimage} to the  wake at $f(c)$ and obtain that $f^i(c)=c$ for some $i>0$ which is impossible  since $c\in\partial B$.
\cqfd

\begin{lem}
For every point $z$ of $\C \setminus B$  there exist a unique $u\in \partial B$  such that $z$ belongs to the  wake $W(u)$.
\end{lem}
\proof This follows directly from the local connectivity or corollary~\ref{c:theta}. \cqfd

\begin{lem}
If $c\in \partial B$, for any $z\in \partial B$ there is exactly one external ray converging to it excepted if there exists $j\ge0$ such that $f^j(z)=c$ in which  case $f^j\from W(z)\to W(c)$ is a  homeomorphism and $W(c)$ contains the preimage of $B$. Moreover every Fatou component is an iterated preimage of $B$.\end{lem}
\proof
This follows directly from Lemma~\ref{l:wakeimage} .\cqfd

\begin{lem}
If $c\notin \partial B$,  for any $z\in \partial B$ there is exactly one external ray converging to it, excepted if there exists $j\ge0$ such that $\Theta(f^j(z))=t$ in which   case $f^j \from W(z) \to W(t)$ is a homeomorphism. 
\end{lem}
\proof It follows from the local connectivity and from  Lemma~\ref{l:wakeimage} \cqfd


\begin{prop}\label{p:periodicrenormalizable}If the map
$f$ has a non-repelling strictly
periodic point $x$  then $f$ is renormalizable
near $x$.
\end{prop}
\proof The sequence $(P_n(x))$ of puzzle pieces containing $x$ is well defined, for any of the graphs considered, since $x$ is not  on the graphs (it is neither eventually repelling and nor $\p$). 
The periodicity of $x$ implies that  there exists $k\ge 1$ such that $f^k(P_{n+k}(x))=P_{n}(x)$ 
 for any large $n$. Take the minimal $k $ satisfying this equality. The critical point $c$ belongs 
  to some piece $P_{n}(f^i(x))$, for every  $n\ge 0$. Otherwise the map $f^k\from P_{n+k}(x)\to P_{n}(x)$ would be invertible and its inverse $g\from P_{n}(x)\to P_{n}(x)$  would have an
attracting fixed point (by Schwarz' Lemma) since 
$ P_{n+k}(x)=g(P_{n}(x))\neq P_{n}(x)$. This is not possible
since  $|g'(x)|\ge 1$. Note that  $i$ is independent of $n$
since the pieces $P_n(f^j(c))$ are
disjoint for $j<k$ (if $k>1$). 
We prove now that if $k>1$, the map $f^k \from
P_{n+k}(f^i(x))\to P_{n}(f^i(x))$ is renormalizable. It is enough to prove that 
$\ol {P_{n+k}(f^i(x))}\subset  P_{n}(f^i(x))$, {\it i.e.} that the critical point 
is surrounded by a non degenerate annulus.
If $c\notin W(0)$, this follows from Lemma~\ref{l:condition}.
If $c\in W(0)$ and if there is $k$ such that $f^k(c)$ is not in $W(0)$ then Lemma~\ref{l:w0} implies the result.
Now, consider the case where  the orbit of $c$ stays in $W(0)$. The rays 
$R(0)$ and $R(1/2)$  bound $W(0)$ and converge to $\p$.
The map $f$ has three fixed points in $\C$, one double at $\p$, so that there is a fixed point $\alpha$ in $W(0)$ (since it contains $c$, see~\cite{GM}). This point is repelling (by hypothesis), so 
there is a cycle of external rays converging to it, say $R(3^j\zeta)$ with $0\le j\le m$.
Consider the connected component  $U$ of $W(0)\setminus (R_\infty(1/6)\cup R_\infty(1/3)\cup R_\infty(2/3)\cup R_\infty(5/6)\cup E_\infty(1))$ containing the critical point.
We define a new puzzle in $U$. Let $\Gamma$ be the graph formed in $U$ by the cycle
$R(3^j\zeta)$ with $0\le j\le m$ together with $\partial U$. Denote by $Y_n$ the puzzle pieces. 
Since $x$ is periodic by $f$, there is some minimal $l$ such that $f^l$ maps $Y_{n+l}(x)$ onto $Y_n(x)$. Moreover this map has degree $2$ otherwise  by Schwarz Lemma the fixed point $x$ would be repelling for $f^l$. Since  $\alpha$ is repelling, one can   thicken  near $\alpha$ the puzzle pieces  $Y_n$  to get larger puzzle pieces $\tilde Y_n$, as done in Milnor's article~\cite{Mi3}, so that $f^l\from \tilde Y_{n+l}(x)\to \tilde Y_n(x)$ is renormalizable. The result follows. \cqfd

\begin{crr}\label{c:brujno} Assume that $f$ possesses a strictly periodic point $p$
with multiplier $\lambda=e^{2i\pi \theta}$, such that $\theta\in
\R\setminus\Q$.
  Then $f$ is linearizable near $p$ if and only if $\theta \in \mathcal B$.
  Moreover, if
$\theta\notin \mathcal B$ there exist periodic cycles in any
neighbourhood of $p$.
\end{crr}
Here $\mathcal B$ denotes the set of {\it Brjuno} numbers\,: an
irrational $\theta$ of convergents $p_n/q_n$ (rational
approximations obtained by the continued fraction development) is
a {\it Brjuno} number,   if $\sum_{n=1}^{\infty}(\log
q_{n+1})/q_n$ is finite.

\proof The map  $f$ is renormalizable by
Proposition~\ref{p:periodicrenormalizable}. So there is a
homeomorphism   that conjugates $f^k$ to a quadratic
polynomial $z^2+c$ on a neighbourhood of the corresponding  Julia
set (see~\cite{DHP}). The multiplier at the fixed points are the
same by Na\"{\i}shul' Theorem (see~\cite{Na}). So the result
follows from Yoccoz' and Brjuno's work (see~\cite{Yo}). \cqfd


\bigskip
\begin{thebibliography}{??}
 \bibitem[Bl]{Blanchard} {\sc P. Blanchard} ---
             {\em  Complex analytic dynamics on the Riemann sphere},
              Bull. Amer. Math. Soc.  {\bf 11} (1984), 85--141.
\bibitem[Ca-Ga]{CG} {\sc L. Carleson, T.W. Gamelin} ---
                {\em Complex Dynamics}, Universitext: Tracts in Mathematics. Springer-Verlag, New York, 1993. x+175 pp. ISBN: 0-387-97942-5 
\bibitem[D-H1]{DH} {\sc A. Douady, J.~H. Hubbard}  ---
             {\em Etude dynamique des polyn\^omes complexes},
             Publications math\'ematiques d'Orsay 1984.
\bibitem[D-H2]{DHP} {\sc A. Douady, J.~H. Hubbard}  ---
             {\em On the dynamics of polynomial-like mappings},
              Ann. scient. \'Ec. Norm. Sup. {\bf 18} (1985), 287--343.

\bibitem[Go-Mi]{GM} {\sc L. Goldberg, J. Milnor}  ---         {\em Fixed points of polynomial maps. II. Fixed point portraits}. Ann. Sci. ƒcole Norm. Sup. (4) 26 (1993), no. 1, 51--98. 58F20 (58F03).

\bibitem[Ki]{Kiwi} {\sc J.~Kiwi} ---
             {\em $\Bbb R$eal laminations and the topological dynamics of complex polynomials.} Adv. Math. 184 (2004), no. 2, 207--267. 
             
 \bibitem[Mi1]{Mi1} {\sc J. Milnor}  ---
            {\em Dynamics in One Complex Variable},
            Vieweg 1999, 2nd edition 2000.
Annals of Mathematics Studies,160. Princeton University Press, Princeton, NJ, 2006.
 \bibitem[Mi2]{Mi2} {\sc J. Milnor}  ---
          {\em On cubic polynomial maps with periodic critical point}, preprint (1991), to
            be published in a volume dedicated to JHH.
\bibitem[Mi3]{Mi3} {\sc J. Milnor}  ---
             {\em Local Connectivity of Julia Sets: Expository Lectures},
              pp. 67-116 of ``The Mandelbrot set, Theme and Variations''
              edit: Tan Lei, LMS Lecture Note Series 274 ,
              Cambr. U. Press 2000.

\bibitem[Na\u\i]{Na} {\sc  V.A. Na\u\i shul'}, ---
             {\em  Topological invariants of analytic and area-preserving
              mappings and their application to analytic differential
             equations in $ C\sp{2}$ and $ CP\sp{2}$. } ---
             (Russian) Funktsional. Anal. i Prilozhen. 14 (1980), no. 1, 73--74. 
\bibitem[Pe1]{Pe1} {\sc C. L. Petersen} ---
            {\em On the Pommerenke-Levin-Yoccoz inequality},
              Ergodic Theory Dynam. Systems 13 (1993), no. 4, 785--806.

\bibitem[Po]{Po} {\sc C. Pommerenke} ---
              {\em Boundary Behaviour of Conformal Maps},
              Springer Verlag 1992.

\bibitem[Ro1]{Ro1} {\sc P.~Roesch} ---
             {\em Puzzles de Yoccoz pour les applications ˆ allure rationnelle},  Enseign. Math. (2) 45 (1999), no. 1-2, 133--168.  
\bibitem[Ro2]{Ro2} {\sc P.~Roesch} ---
             {\em Hyperbolic components 
             of polynomials with fixed critical point of maximal
             order},  Annales de L'ENS, to appear.
\bibitem[Ro3]{Ro3} {\sc P.~Roesch} ---
             {\em  Cubic parabolic slice},  manuscript.
            
\bibitem[TaYi]{TY} {\sc Tan Lei ,  Yin Yongcheng } ---
              {\em Local connectivity of the Julia set for geometrically finite
              rational maps}, Science in China (Serie A) {\bf 39} (1996),
              39--47.

\bibitem[Yo]{Yo} {\sc J.-C. Yoccoz} ---
             {\em Petits diviseurs en dimension $1$}, Ast\'erisque {\bf 231}, 1995.

\end{thebibliography}
\end{document}